\documentclass[3p,twocolumn]{elsarticle}

\usepackage{graphics} 
\usepackage{epsfig} 
\usepackage{amsmath} 
\usepackage{amssymb}  
\usepackage{cite}
\usepackage{color}
\usepackage{graphicx}
\usepackage{algpseudocode}
\usepackage{algorithm}
\usepackage{amsthm}

\usepackage{enumitem}

\usepackage[font=footnotesize]{caption}

\newtheorem{theorem}{Theorem}[section]
\newtheorem{definition}[theorem]{Definition}

\newtheorem{lemma}[theorem]{Lemma}

\newtheorem{remark}[theorem]{Remark}

\newtheorem{corollary}[theorem]{Corollary}
\newtheorem{proposition}[theorem]{Proposition}  
\newtheorem{problem}[theorem]{Problem}

\newcommand{\naturals}{\mathbb{N}}
\newcommand{\real}{\mathbb{R}}
\newcommand{\cplx}{\mathbb{C}}

\newcommand{\Fc}{\mathcal{F}}

\newcommand{\Ic}{\mathcal{I}}
\newcommand{\Jc}{\mathcal{J}}
\newcommand{\Kc}{\mathcal{K}}
\newcommand{\Lc}{\mathcal{L}}

\newcommand{\Pc}{\mathcal{P}}

\newcommand{\Sc}{\mathcal{S}}

\newcommand{\Vc}{\mathcal{V}}
\newcommand{\Wc}{\mathcal{W}}
\newcommand{\Xc}{\mathcal{X}}

\newcommand{\Fbb}{\mathbb{F}}

\newcommand{\Kapprox}{K_{\operatorname{approx}}}
\newcommand{\apprx}{\operatorname{approx}}

\newcommand{\until}[1]{\{1,\dots,#1\}}

\newcommand{\Span}{\operatorname{span}}

\newcommand{\restr}[2]{#1 \!\! \restriction_{#2}}

\newcommand{\innerprod}[2]{\langle #1, #2 \rangle}

\newcommand{\longthmtitle}[1]{\mbox{}{\textit{(#1):}}}
\newcommand{\setdef}[2]{\{#1 \; | \; #2\}}

\newcommand{\oprocendsymbol}{\hbox{$\square$}}
\newcommand{\oprocend}{\relax\ifmmode\else\unskip\hfill\fi\oprocendsymbol}
\def\eqoprocend{\tag*{$\square$}}

\allowdisplaybreaks

\graphicspath{{epsfiles/}}
\begin{document}

\begin{frontmatter}

\title{\bf Invariance Proximity: Closed-Form Error Bounds
for Finite-Dimensional Koopman-Based Models}

\author[1]{Masih Haseli\corref{cor1}}
\ead{mhaseli@ucsd.edu}
\author[1]{Jorge Cort{\'e}s}
\ead{cortes@ucsd.edu}
\cortext[cor1]{Corresponding author}
\address[1]{Department of Mechanical and Aerospace Engineering,
University of California, San Diego, 9500 Gilman Dr, La Jolla, CA 92093}


\begin{abstract}
A popular way to approximate the Koopman operator's action on a
finite-dimensional subspace of functions is via orthogonal
projections. The quality of the projected model directly depends on
the selected subspace, specifically on how close it is to being
invariant under the Koopman operator. The notion of invariance
proximity provides a tight upper bound on the worst-case relative
prediction error of the finite-dimensional model. However, its
direct calculation is computationally challenging.  This paper
leverages the geometric structure behind the definition of
invariance proximity to provide a closed-form expression in terms of
Jordan principal angles on general inner product spaces.  Unveiling
this connection allows us to exploit specific isomorphisms to
circumvent the computational challenges associated with spaces of
functions and enables the use of existing efficient numerical
routines to compute invariance proximity.
\end{abstract}
\end{frontmatter}

\section{Introduction}
Koopman operator theory studies the behavior of nonlinear dynamical
systems through the lens of linear operators acting on vector spaces
of functions. This paradigm provides a formal algebraic structure that
can be leveraged to study unstructured complex systems. However, the
Koopman operator is generally defined on infinite-dimensional spaces,
a major obstruction for implementation on digital computers. A popular
way to address this is to approximate the action of the operator over
finite-dimensional subspaces. Expectedly, such approximations, often
calculated via orthogonal projections, lead to model mismatch and
prediction error, which makes providing accuracy measures for such
models critical. These measures can be employed both as loss functions
for subspace learning and as a tool to provide error bounds and
certificates on accuracy and safety for the actual system. This paper
focuses on the computation of one such measure: invariance proximity.

\subsection{Literature Review}
Koopman operator theory~\citep{BOK:31} describes a nonlinear system via
the action of a linear operator on a vector space of
functions. Moreover, the value of Koopman eigenfunctions on the system
trajectories evolve linearly. This leads to a powerful spectral
representation of nonlinear systems~\citep{IM:05}, which has given rise
to a plethora of applications, including stability
analysis~\citep{AM-IM:16,JJB-GF:24,SAD-AMV-CJT:22},
control~\citep{MK-IM-automatica:18,SP-SK:19,VZ-EB:23-neurocomputing,MH-JC:24-auto,DU-KD:23,RS-JB-FA:23},
and robotics~\citep{DB-XF-RBG-CDR-RC:20,GM-IA-TDM:23}.
In inner-product spaces, one can approximate the action of the Koopman
operator on finite-dimensional subspaces by using orthogonal
projections.  The description of the evolution of observables
under the operator naturally lends itself to the incorporation of data
in producing such approximations.  A notable example is Extended
Dynamic Mode Decomposition (EDMD)~\citep{MOW-IGK-CWR:15}, which uses
data to approximate the Koopman operator's action on a predefined
finite-dimensional space spanned by a dictionary of functions.  
The work~\citep{MK-IM:18} provides several convergence results
regarding EDMD's behavior with respect to the action of the Koopman
operator as the number of data and the dictionary functions go to
infinity. Notably, for bounded Koopman operators, these results
include convergence of EDMD operators in strong operator topology,
weak convergence of eigenfunctions, as well as important
implications for \emph{finite horizon} accurate prediction of
observables on system's trajectories. We refer the reader
to~\citep{MJK:23} for a survey on different DMD-based methods and
their properties.

In practical applications, due to computational constraints, one
often relies on approximations on finite-dimensional spaces.  The
accuracy of EDMD and other projection-based methods directly depends
on the quality of the underlying finite-dimensional subspace. If the
subspace is invariant under the Koopman operator, the resulting
model is exact. Otherwise, the approximation via projection leads to
two related issues: (i) errors in the prediction of the operator's
action and (ii) spectral pollution. This has led to significant
research activity towards finding finite-dimensional spaces that are
(close to) invariant under the action of the Koopman operator. The
works~\citep{NT-YK-TY:17,BL-JNK-SLB:18,MS:21-L4DC} use optimization
and neural network-based methods to address this question. Other
methods directly identify the spectral properties of the Koopman
operator, including its eigenfunctions, which in turn span invariant
subspaces~\citep{MK-IM:20,EK-JNK-SLB:21}. In another line of
work~\citep{MH-JC:22-tac,MH-JC:21-tcns}, termed Symmetric Subspace
Decomposition (SSD) algorithms, we have relied on iteratively
cleaning up a given subspace to identify the maximal Koopman
invariant subspace and all the eigenfunctions with convergence and
accuracy guarantees. The Tunable Symmetric Subspace Decomposition
(T-SSD) algorithm~\citep{MH-JC:21-acc,MH-JC:23-auto} allows to
perform the clean-up with tunable accuracy. The accuracy in T-SSD is
captured via the concept of invariance proximity. T-SSD also gives
convergence and accuracy guarantees on the prediction of the Koopman
operator's action for \emph{all} functions (not just the
eigenfunctions) in the identified subspace.

Residual Dynamic Mode Decomposition
(ResDMD)~\citep{MJC-LJA-MS:23,MJC-AT:24,MJK:24} is a more recent
line of work which aims at resolving the issue of spectral
pollution. ResDMD relies on a clean-up procedure based on the
concept of residuals to remove spurious eigenfunctions and also
provides methods for approximation of Koopman operators'
psudospectra accompanied by appropriate spectral bounds.  Albeit
ResDMD, SSD, and T-SSD all employ some form of clean-up procedure,
it is worth noting that the goals for which they are designed and
the associated guarantees are different.  The clean-up in ResDMD is
focused on the accuracy of candidate eigenpairs, since the goal is
to capture the spectrum. In SSD and T-SSD algorithms, the clean-up
procedure is focused on ensuring the accuracy of prediction for
Koopman operator's action on \emph{arbitrary} functions (not just
eigenfunctions). Moreover, the type of guarantees for ResDMD, SSD,
and T-SSD are different. ResDMD leads to spectral bounds and
guarantees on capturing spectra and psuedospectra, while SSD and
T-SSD lead to bounds on the approximation of the Koopman operator's
action on arbitrary functions in their identified subspace. Another
notable work for spectral bounds is~\citep{VK-KL-PN-MP:23} which
studies systems with compact and self-adjoint Koopman operators.

Given the recent interest in finite-dimensional Koopman-based
approximations, directly characterizing the approximation's accuracy
of given model is critical for the validation and refinement of
Koopman-based models.
The work in~\citep{FN-SP-FP-MS-KW:23} provides probabilistic error
bounds for accuracy of EDMD based on sampled data
and~\citep{MH-JC:23-csl} provides a tight upper-bound for the error
induced by EDMD's projection.
However, these bounds are only given for data-driven techniques, and
not in the larger context of Koopman operator-based methods.
Nonetheless, in system and control theory, many of the applications
(e.g.,~stability, reachability, safety analysis, identification of
invariant sets) require analytical bounds at the function level, not
just on the data.
We tackle this by providing error bounds on Koopman-based projected
models over general inner product spaces.

\subsection{Statement of Contributions}
Given a general inner-product space,
we consider approximate Koopman-based models created by the orthogonal
projection on finite-dimensional subspaces.  To assess the accuracy of
such models, we rely on the notion of invariance proximity introduced
in~\citep{MH-JC:24-auto}, which measures the worst-case relative error
for the model's prediction over all functions in the
subspace. Given that invariance proximity requires taking a
supremum over uncountably many functions in a finite-dimensional
vector space, its \emph{direct} calculation is challenging.
To efficiently compute invariance proximity\footnote{In this paper, we
consider the problem of efficiently computing invariance
proximity. We refer the reader to~\citep{MH-JC:24-auto} for the
theoretical properties of invariance proximity and how it can be
used to approximate Koopman-based models.}, we study the geometric
structure between the underlying subspace and its image under the
Koopman operator. We orthogonally decompose these subspaces using
their corresponding Jordan principal angles and vectors. Using this
decomposition, we provide a closed-form expression for invariance
proximity as the sine of the largest principal angle.  This allows us
to exploit specific isomorphisms to reformulate on a complex Euclidean
space the problem of calculating invariance proximity, enabling the
use of existing efficient numerical routines.

\subsection{Notations}
We use $\naturals$, $\real$, $\real_{\geq 0}$, and $\cplx$ to
represent natural, real, non-negative real, and complex numbers. Given
the vector $v \in \cplx^ n$, we denote its complex conjugate, norm,
and transpose with $\bar{v}$, $\|v\|$, and $v^T$ respectively. Given
sets $A$ and $B$, $A \subseteq (\subset) B$ means that $A$ is a
(proper) subset of $B$. Given the vector space $\Lc$ and subspaces
$\Vc,\Wc \subseteq \Lc$, we define their sum
$\Vc+\Wc := \setdef{v + w}{v \in \Vc, \, w \in \Wc}$. Moreover, if
$\Lc$ is equipped with an inner product, $\Vc \perp \Wc$ means that
$\Vc$ is orthogonal to $\Wc$. In this case, we denote their sum with
$\Vc \oplus \Wc$ and refer to it as a direct sum. Given,
$\theta \in \real$, we show its sine and cosine by $\sin(\theta)$ and
$\cos(\theta)$.  Given the functions $f_1$ and $f_2$ with matched
domains and co-domains, we denote their composition by
$f_1 \circ f_2$.

\section{Preliminaries}\label{sec:preliminaries}
We briefly recall~\citep{MB-RM-IM:12,MH-JC:24-auto} the definition of
the Koopman operator, finite-dimensional approximations, and accuracy
bounds characterized through invariance proximity.

\subsection{Koopman Operator}\label{sec:prelim-Koopman}
Consider a discrete-time nonlinear system
\begin{align}\label{eq:dynamical-sys}
x^+ = T(x), \; x \in \Xc,
\end{align}
where $\Xc$ is the state space and $T: \Xc \to \Xc$ is the dynamics
map. Consider a vector space $\Fc$ over $\cplx$ comprised of
complex-valued functions with domain $\Xc$ and assume it is closed under
composition with $T$: for all $f \in \Fc$, we have
$f \circ T \in \Fc$. The Koopman operator $\Kc: \Fc \to \Fc$ is
defined as
\begin{align}\label{eq:Koopman-def}
\Kc f = f \circ T.
\end{align}
Unlike the system~\eqref{eq:dynamical-sys}, which acts on points in
the state space, the Koopman operator~\eqref{eq:Koopman-def} acts on
functions in the vector space $\Fc$.
Importantly, the Koopman operator is always linear, i.e., for all
$g,h \in \Fc$ and all $\alpha, \beta \in \cplx$,
\begin{align}\label{eq:Koopman-spatial-linear}
\Kc (\alpha g + \beta h) = \alpha \,\Kc g + \beta \, \Kc h.
\end{align}
A nonzero function $\phi \in \Fc$ is an eigenfunction of the Koopman
operator with eigenvalue $\lambda \in \cplx$ if
\begin{align}\label{eq:Koopman-eig-def}
\Kc \phi = \lambda \phi.
\end{align}
The eigenfunctions evolve linearly in time on the system's
trajectories, i.e.,
$ \phi(x^+) = \phi \circ T(x) = [\Kc \phi] (x) = \lambda \phi (x)$.
The combination of the linear temporal evolution of eigenfunctions
with the linearity~\eqref{eq:Koopman-spatial-linear} of the operator
lead to computationally efficient methods for identification and
prediction of nonlinear systems.  It is crucial to note that, in
general, the space $\Fc$ is infinite-dimensional.  For many practical
application, a finite-dimensional representation is used, as we
explain next.

\subsection{Koopman-Invariant Subspaces}
A subspace $\Jc \subset \Fc$ is invariant under the Koopman operator
if $\Kc f \in \Jc$ for all $f \in \Jc$.  Although
finite-dimensional Koopman-invariant subspaces capturing complete
information about the system are rare in general, their study is
theoretically important because the concept of subspace invariance
determines the general form of finite-dimensional models and
provides a bedrock for approximations on non-invariant
subspaces. Finite-dimensional Koopman-invariant subspaces
allow for
exact representation of the Koopman operator and enable the use of
efficient numerical linear algebraic routines. This exact
representation is constructed by restricting $\Kc$ to a
finite-dimensional subspace $\Jc$ as
$ \restr{\Kc}{\Jc}: \Jc \to \Jc $, where $\restr{\Kc}{\Jc} f = \Kc f$
for all $f \in \Jc$.  Given a basis for $\Jc$, one can represent the
operator $\restr{\Kc}{\Jc}$ by a matrix. Formally, let
$J: \Xc \to \cplx^{\dim{\Jc}}$ be a vector-valued map whose elements
form a basis for $\Jc$. Then, there exists
$K \in \cplx^{\dim{\Jc} \times \dim{\Jc}}$ such that
\begin{align}\label{eq:restriction-on-invariant-basis}
\restr{\Kc}{\Jc} J = \Kc J = J \circ T = K J,
\end{align}
where the action of an operator on a vector-valued map is defined in
an element-wise manner.  For any function $f \in \Jc$ represented as
$f = v_f^T J$, the action of $\restr{\Kc}{\Jc}$ on $f$ is
\begin{align}\label{eq:restriction-on-invariant-function}
\restr{\Kc}{\Jc} f = v_f^T K J.
\end{align}
Equations~\eqref{eq:restriction-on-invariant-basis}-\eqref{eq:restriction-on-invariant-function}
enable fast prediction of the action of the Koopman operator via
numerical linear algebra. However, in general, finding
finite-dimensional Koopman-invariant subspaces that capture sufficient
information is difficult (and sometimes impossible) and therefore, one
settles for approximations, as we explain next.

\subsection{Approximations on Non-Invariant Subspaces}
To discuss approximations to the Koopman operator on non-invariant
subspaces, throughout the paper we equip the space $\Fc$ with an inner
product
$\innerprod{\cdot}{\cdot}: \Fc \times \Fc \to \cplx$ which induces the
norm $\| \cdot \|: \Fc \to \real_{\geq 0}$.  In this paper, we do
\emph{not} assume the space $\Fc$ is complete (Hilbert) to allow for
more general settings. We aim to approximate the action of the
Koopman operator on a finite-dimensional space $\Sc \subset \Fc$ which
is not Koopman invariant. To tackle this, consider the orthogonal
projection operator $\Pc_{\Sc}: \Fc \to \Fc$ on $\Sc$, which maps a
function in $\Fc$ to the closest function in $\Sc$.  Before
proceeding, we first remark that even though $\Fc$ might not be
complete (Hilbert), the best approximation on the finite-dimensional
subspace $\Sc$ always exists and is unique, therefore the orthogonal
projection operator $\Pc_{\Sc}$ is well defined.

\begin{remark}\longthmtitle{Existence and Uniqueness of Best
Approximations on Finite-dimensional Subspaces} {\rm
Finite-dimensional subspaces of an inner product space (on
fields $\real$ or $\cplx$) are Chebyshev sets, that is, every
point in the inner-product space has a unique closest point on
the finite-dimensional subspace, see
e.g.,~\citep[Result~3.8~(3)]{FD:01}.
This is a direct consequence of the completeness of
finite-dimensional subspaces of inner-product spaces (which
might not be complete themselves) on fields $\real$ or $\cplx$,
which can be derived from~\citep[Theorem~2.4-2]{EK:89}. In fact,
the closest point coincides with the orthogonal projection, which
can be computed in closed form: given the finite-dimensional
subspace $\Sc \subseteq \Fc$ and an arbitrary function
$f \in \Fc$, the orthogonal projection of $f$ on $\Sc$ can be
computed as $\Pc_{\Sc} f = \sum_{i=1}^n \innerprod{f}{e_i}e_i$,
where $\{e_1,\ldots,e_n\}$ is an \emph{orthonormal} basis for
$\Sc$.  \oprocend }
\end{remark}

To approximate the action of the Koopman operator on the
finite-dimensional subspace $\Sc$ in a way that allows to work with
matrix representations, we approximate the Koopman operator $\Kc$
by $\Kc_{\apprx} := \Pc_{\Sc} \Kc: \Fc \to \Fc$. Note that the space
$\Sc$ is invariant under $\Kc_{\apprx}$; hence, given a basis for
$\Sc$ represented by the vector-valued map
$\Psi: \Xc \to \cplx^{\dim(\Sc)}$, we can
apply~\eqref{eq:restriction-on-invariant-basis} to $\Kc_{\apprx}$ (by
swapping $\Kc$, $\Jc$, and $J$ with $\Kc_{\apprx}$, $\Sc$, and $\Psi$
resp.) to approximate the action of the Koopman operator as:
\begin{align}\label{eq:restriction-on-noninvariant-basis}
\Kc \Psi = \Psi \circ T \approx \restr{\Kc_{\apprx}}{\Sc} \Psi
=  \Kapprox \Psi, 
\end{align}
where $\Kapprox \in \cplx^{\dim(\Sc) \times \dim(\Sc)}$. Moreover,
similarly to~\eqref{eq:restriction-on-invariant-function}, for any
function $f \in \Sc$ with representation $f = v_f^T \Psi$, one can
approximate the action of the Koopman operator on $f$ by
\begin{align}\label{eq:restriction-on-noninvariant-function}
\Kc f \approx \Kc_{\apprx} f =	\restr{\Kc_{\apprx}}{\Sc}  f =
v_f^T \Kapprox \Psi. 
\end{align}

\begin{remark}\longthmtitle{Connections to Extended Dynamic Mode
Decomposition
(EDMD)~\citep{MOW-IGK-CWR:15}}\label{r:EDMD-projection}
{\rm The EDMD method is a special case of approximations
in~\eqref{eq:restriction-on-noninvariant-basis}-\eqref{eq:restriction-on-noninvariant-function},
where $\Psi$ is the selected dictionary and $\Fc$ is the space
$L^2$ defined with respect to the empirical measure on the data
set, see e.g.,~\citep{MK-IM:18,SK-PK-CS:16}.  \oprocend }
\end{remark}

The quality of approximation
in~\eqref{eq:restriction-on-noninvariant-basis}-\eqref{eq:restriction-on-noninvariant-function}
directly depends on the quality of the subspace $\Sc$. If $\Sc$ is
invariant under the operator,
equations~\eqref{eq:restriction-on-noninvariant-basis}-\eqref{eq:restriction-on-noninvariant-function}
reduce
to~\eqref{eq:restriction-on-invariant-basis}-\eqref{eq:restriction-on-invariant-function}
and there is no approximation error. Otherwise, the projection
in~\eqref{eq:restriction-on-noninvariant-basis}-\eqref{eq:restriction-on-noninvariant-function}
leads to information loss and approximation error. In practical
applications, one requires bounds on the quality of the model;
therefore, it is of utmost importance to quantify the approximation
accuracy, as we discuss next.

\subsection{Invariance Proximity}
Here, we present the concept of invariance proximity
following~\citep[Definition~8.5]{MH-JC:24-auto} to measure the quality
of a finite-dimensional subspace $S \subset \Fc$ in terms of how close
it is to being invariant under the Koopman operator. Invariance
proximity is formally given by
\begin{align}\label{def:invariance-proximity}
\Ic_{\Kc} (\Sc)
&:= \sup_{f \in \Sc, \| \Kc f \| \neq 0} \frac{\|
\Kc f - \Pc_{\Sc} \Kc f\|}{ \| \Kc f \|}  \notag 
\\
&=  \sup_{f \in \Sc, \| \Kc f \| \neq 0} \frac{\| \Kc f -
\Kc_{\apprx} f\|}{ \| \Kc f \|} .
\end{align}
This measures the worst-case relative error of the
approximation~\eqref{eq:restriction-on-noninvariant-function} of the
operator's action.
It only depends on the operator $\Kc$ and the subspace $\Sc$ (since
$\Kc_{\apprx} = \Pc_{\Sc} \Kc$ only depends on $\Kc$ and $\Sc$), and
\emph{does not} depend on the choice of basis for~$\Sc$.  Even though
out of scope of this paper, we note that invariance proximity also
provides bounds on Koopman-based models for control systems,
cf.~\citep{MH-JC:24-auto}.
Despite its importance, there do not exist methods to compute
invariance proximity in general inner product spaces. For the
particular case of 
the space $L^2$ with respect to the empirical measure on a data set,
where the EDMD method operates, cf.~Remark~\ref{r:EDMD-projection},
one can obtain~\citep{MH-JC:23-csl} a closed-form expression for
invariance proximity based on the application EDMD forward and
backward in time.
Given the important properties of invariance proximity~in~general
inner product spaces, the problem considered here is the~development
of efficient methods to compute it in a general setting.

\subsection{Bounding Eigenspace Residuals via Invariance
Proximity} On a first look, the ratio in the definition of
invariance proximity might look similar to the definition of
residuals in ResDMD, cf.~\citep[Equations~(3.1) and
(3.4)]{MJC-LJA-MS:23}. However, these notions turn out to be
significantly different. Specifically, with the notation adopted
here, for a candidate eigenpair $(\lambda, \phi)$ of the Koopman
operator, the residual of this pair is defined as\footnote{The
residuals in~\citep{MJC-LJA-MS:23} are defined on space $L^2$;
however, they can easily be extended to general inner-product
spaces as long as the operators involved are well defined.}
\begin{align*}
\operatorname{res}(\lambda,\phi) = \frac{\| \Kc \phi - \lambda \phi
\|}{\| \phi \|}. 
\end{align*}
Note that both the numerator and the denominator of the residuals are
different from the ratios in the definition of invariance
proximity. The numerator in~\eqref{def:invariance-proximity} is
defined for \emph{arbitrary} functions, while residuals are only
defined for a candidate \emph{eigenpair}. Even if we restrict the
vector space $\Sc$ in invariance
proximity~\eqref{def:invariance-proximity} to a one-dimensional space
spanned by an eigenfunction of $\Kc_{\apprx}$, the notions are still
different, as we explain next.

If the eigenpair is chosen based on the operator
$\restr{\Kc_{\apprx}}{\Sc}: \Sc \to \Sc$
in~\eqref{eq:restriction-on-noninvariant-basis}, i.e., if
$\restr{\Kc_{\apprx}}{\Sc} \phi = \lambda \phi$, then the residual
above turns into
\begin{align}\label{eq:residul}
\operatorname{res}(\lambda,\phi) = \frac{\| \Kc \phi - \Kc_{\apprx}
\phi \|}{\| \phi \|}. 
\end{align}
Note that~\eqref{eq:residul} is \emph{not} a special case of
invariance proximity for the one-dimensional subspace $\Span(\phi)$,
since the denominator of~\eqref{def:invariance-proximity} is
$\|\Kc f\|$ not $\|f\|$. It is worth mentioning that the residual
in~\eqref{eq:residul} is not a relative error for prediction of
$\Kc \phi$ and in a general setting its value can exceed one, while
invariance proximity is always between zero and one, see
e.g.,~\citep[Proposition~8.6]{MH-JC:24-auto}.

Invariance proximity and residuals have different uses: invariance
proximity measures the accuracy of the model built based on the
Koopman operator on a finite-dimensional space and leads to explicit
tight bounds on prediction of Koopman operator's \emph{action} on all
(uncountably many) functions in the space, see,
e.g.,~\citep{MH-JC:21-acc,MH-JC:23-auto} for its early use for
subspace identification in the context of dynamic mode decomposition
(the special case of space $L^2$ on empirical measure,
cf.~Remark~\ref{r:EDMD-projection}).  Residuals on the other hand,
provide a way to approximate the spectra and pseudospectra of the
Koopman operator.  The problems of prediction of the operator's action
on arbitrary functions and approximation of spectra are different and
hence require different tools to handle.

Even though invariance proximity and residuals are not related in a
general setting, for the special case where the Koopman operator is
bounded, invariance proximity can provide a bound for the residuals
in~\eqref{eq:residul}.

\begin{lemma}\longthmtitle{Bounding Residuals via Invariance Proximity}
Let the Koopman operator $\Kc$ in~\eqref{eq:Koopman-def} be
bounded. Moreover, consider the finite-dimensional subspace $\Sc$,
and the operator $\restr{\Kc_{\apprx}}{\Sc}: \Sc \to \Sc$
in~\eqref{eq:restriction-on-noninvariant-basis}. Then, for any
eigenpair $(\lambda,\phi)$ of $\restr{\Kc_{\apprx}}{\Sc}$ satisfying
$\restr{\Kc_{\apprx}}{\Sc} \phi = \lambda \phi$ and
$\| \Kc \phi \| \neq 0$, we have
\begin{align*}
\operatorname{res}(\lambda,\phi)  \leq \| \Kc \| \, \Ic_{\Kc}(\Sc).
\end{align*}
\end{lemma}
\begin{proof}
By the definition of invariance
proximity~\eqref{def:invariance-proximity}, one can write
\begin{align*}
\frac{\| \Kc \phi - \Kc_{\apprx} \phi \|}{\| \Kc \phi \|} \leq
\Ic_{\Kc}(\Sc). 
\end{align*}
Moreover, since the Koopman operator is bounded we have
$\| \Kc \phi \| \leq \| \Kc \| \| \phi \|$, which in conjunction
with the equation above and~\eqref{eq:residul} leads to
\begin{align*}
\operatorname{res}(\lambda,\phi)
&= \frac{\| \Kc \phi - \Kc_{\apprx} \phi \|}{\| \phi \|} 
\\
&\leq  \| \Kc\| \frac{\| \Kc \phi - \Kc_{\apprx} \phi \|}{\| \Kc
\phi \|}  \leq  \| \Kc \| \, \Ic_{\Kc}(\Sc), 
\end{align*}
which concludes the proof.
\end{proof}

\section{Problem Statement}
The importance of invariance proximity stems from the fact that it
provides a tight upper bound on the worst-case error induced by
projecting the action of the Koopman operator on a finite-dimensional
space. Therefore, it can be used to find error bounds of Koopman-based
models for both data-driven and analytic prediction and control of
dynamical processes. In addition, one can use invariance proximity as
a cost function for optimization-based subspace learning, enabling the
identification of Koopman eigenfunctions and eigenmodes.  All the
aforementioned applications are only possible if one can compute
invariance proximity in the inner-product space where the projection
of the action of the Koopman operator is done. Our aim is to
efficiently compute invariance proximity.

\begin{problem}\longthmtitle{Computing Invariance Proximity}
Given the system~\eqref{eq:dynamical-sys} and its associated Koopman
operator~\eqref{eq:Koopman-def} on the \emph{general} inner-product
space $\Fc$, efficiently compute invariance proximity,
$\Ic_{\Kc} (\Sc)$, for any \emph{arbitrary finite-dimensional
subspace} $\Sc \subset \Fc$.  \oprocend
\end{problem}

\section{Principal Angles and Vectors}\label{sec:prinicipal-angles}

Our starting point to unveil the geometric and algebraic structures
behind the notion of invariance proximity is the observation that,
from~\eqref{def:invariance-proximity}, one can see that it depends on
the projection onto the finite-dimensional space $\Sc$ of the image of
$\Sc$ under the Koopman operator~$\Kc$. Therefore, we need tools to
study the relation between two subspaces with respect to each other
and provide algebraic decompositions which simplify working with
projections. To do so, we rely on the well-known notion of principal
angles between vector spaces, initially defined by
Jordan~\citep{CJ:1875} and later formalized in~\citep{HH:92}.  Our
problem setting here is slightly different from existing notions in
the literature because we are interested in finite-dimensional
subspaces of an infinite-dimensional complex inner product space,
which is not necessarily Hilbert.  Given this difference, we cannot
directly use the existing results, and
therefore provide the necessary definitions and results for our
problem setting. In our exposition, we keep the terminology close
to~\citep{AG-CSJH:06}, which studies principal angles in
finite-dimensional complex spaces\footnote{One could indirectly use
the structure in~\citep{AG-CSJH:06} via (uncountably many)
isomorphisms: however, we avoid this route for ease of exposition.}.

\begin{definition}\longthmtitle{Principal Angles Between Subspaces}
\label{def:principal-angles}
Consider\footnote{All results in the paper are valid if $\Fc$ is
defined over the field $\real$ provided that one replaces the
inner product with a real-valued inner product.} two
finite-dimensional subspaces $U, V \subseteq \Fc$ and, without loss
of generality, assume $m_1 := \dim(U) \geq \dim(V) =: m_2$.
Then, the \emph{(Jordan) principal angles}
$0 \leq \theta_1 \leq \cdots \leq \theta_{m_2} \leq \frac{\pi}{2}$
and their corresponding \emph{principal unit vectors},
$\{u_i\}_{i=1}^{m_2} \subset U$ and $\{v_i\}_{i=1}^{m_2} \in V$ are
defined iteratively as:
\begin{align}\label{eq:recursive-principal-angles}
\cos(\theta_i) :=
&\max_{u \in U} \max_{v \in V} |\innerprod{u}{v} | =: \innerprod{u_i}{v_i}
\\
\text{subject to:} \;
&\innerprod{u}{u_k} = 0, \innerprod{v}{v_k} =0, \; \forall k \in \until{i-1}
\nonumber
\\
&\|u\| = 1, \|v\| = 1. \eqoprocend
\end{align}
\end{definition}
\vspace*{-.5ex}

Principal angles depend on the maximum value of the cost function
in~\eqref{eq:recursive-principal-angles} and therefore are
unique. However, principal vectors depend on the maximizers and are
not unique (e.g., $u$ can be replaced with $-u$).

\begin{remark}\longthmtitle{Rotation of Principal Vectors} {\rm In
Definition~\ref{def:principal-angles}, the first set of
constraints is empty for $i=1$. Note that $u_i$s and $v_i$s are
chosen such that their inner product is real, despite being
defined on the field of complex numbers.  Such vectors always
exists since, if $u$ and $v$ are a solution for the optimization
problem~\eqref{eq:recursive-principal-angles}, then one can rotate
$u$ by multiplying with $e^{j \gamma}$ given an appropriate
$\gamma$ to ensure
$\innerprod{e^{j \gamma} u}{v} = |\innerprod{u}{v}|$. Since
$| e^{j \gamma} | = 1$, we can choose $e^{j \gamma} u$ as $u_i$
and $v$ as $v_i$.  \oprocend }
\end{remark}

We state several useful properties of principal vectors.

\begin{lemma}\longthmtitle{Orthonormality of Principal
Vectors}\label{l:orthonormality-principal-vectors}
Given Definition~\ref{def:principal-angles}, for all
$i,j \in \until{m_2}$, we have $\innerprod{u_i}{u_j} = \delta_{ij}$
and $\innerprod{v_i}{v_j} = \delta_{ij}$, where $\delta_{ij}$ is the
Kronecker delta.  \oprocend
\end{lemma}
\smallskip

The proof of this result trivially follows from the constraints
in~\eqref{eq:recursive-principal-angles}.  As a consequence of
Lemma~\ref{l:orthonormality-principal-vectors},
$\{v_1, \ldots, v_{m_2} \}$ is an orthonormal basis for~$V$. Moreover,
if $m_1 = m_2$, then $\{u_1, \ldots, u_{m_2} \}$ is an orthonormal
basis for~$U$.  In case $m_1 > m_2$, we can always add vectors
$\{u_{m_2 +1}, \ldots, u_{m_1}\}$ such that
$\{u_1, \ldots, u_{m_1} \}$ is an orthonormal basis for~$U$.  We use
this convention throughout the paper.

Next, we state an important property regarding the inner product of
principal vectors and the corresponding angles.

\begin{proposition}\longthmtitle{Inner Product of Principal Vectors of
Subspaces}\label{p:innerprod-principalvectors}
The principal vectors satisfy
$\innerprod{u_i}{v_j} = \delta_{ij} \cos(\theta_i)$ for all
$i \in \until{m_1}$ and $j \in \until{m_2}$.
\end{proposition}
\begin{proof}
For the case $i=j \in \until{m_2}$, the equality
$\innerprod{u_i}{v_j} = \cos(\theta_i)$ is a direct consequence of
Definition~\ref{def:principal-angles}. Hence, we only need to prove
$\innerprod{u_i}{v_j} = 0$ when $i \neq j$. We do this by
contradiction. Suppose that $\innerprod{u_i}{v_j} = c \neq 0$.
Consider
\begin{align*}
p_j = \frac{u_j + \bar{c} u_i}{\sqrt{1 + |c|^2}} \in U.
\end{align*}
Using Lemma~\ref{l:orthonormality-principal-vectors}, note that
$\|p_j\| = 1$. The vectors $p_j \in U$ and $v_j \in V$ satisfy the
constraints in the $j$th step of
Definition~\ref{def:principal-angles} (by replacing $i$ with $j$ in
the cost function and constraints
of~\eqref{eq:recursive-principal-angles}). Moreover, one can
write
\begin{align}\label{eq:inner-p-v}
\innerprod{p_j}{v_j} =\frac{1}{\sqrt{1 + |c|^2}} \innerprod{u_j +
\bar{c} u_i}{v_j} = \frac{\innerprod{u_j}{v_j} + |c|^2}{\sqrt{1 +
|c|^2}}, 
\end{align}
where in the last equality we have used the fact that
$\innerprod{u_i}{v_j} = c $. To reach our desired contradiction, we
prove that
\begin{align}\label{eq:contradiction-innerprod}
\innerprod{p_j}{v_j} > \innerprod{u_j}{v_j}.
\end{align}
The previous inequality trivially holds if $\innerprod{u_j}{v_j} = 0$.
Suppose instead that $\innerprod{u_j}{v_j} = r \neq 0$ and note, by
Definition~\ref{def:principal-angles}, that
\begin{align*}
r = \innerprod{u_j}{v_j} =  \cos(\theta_j) \in (0,1].
\end{align*}
Hence, $2r - r^2 >0$ and, consequently, $|c|^4 + (2r-r^2)|c|^2 >
0$. Therefore, adding $r^2 + r^2 |c|^2$ to both sides of the
inequality,
\begin{align}
& |c|^4 + 2 r |c|^2 + r^2 > r^2 + r^2 |c|^2
\nonumber 
\\
&\Rightarrow (|c|^2 + r)^2 > r^2 (1+ |c|^2)
\nonumber
\\ 
&\Rightarrow |c|^2 + r
> r \sqrt{1 + |c|^2} 
\Rightarrow \frac{r + |c|^2}{\sqrt{1 + |c|^2}} > r.
\end{align}
This inequality, in conjunction with $r = \innerprod{u_j}{v_j}$ and
equation~\eqref{eq:inner-p-v}, lead
to~\eqref{eq:contradiction-innerprod}. However,~\eqref{eq:contradiction-innerprod}
directly contradicts the fact $\innerprod{u_j}{v_j}$ is the maximum
in optimization~\eqref{eq:recursive-principal-angles}. Therefore,
the assumption $\innerprod{u_i}{v_j} = c \neq 0$ is false and
$\innerprod{u_i}{v_j} = 0$, which completes the proof.
\end{proof}

An important consequence of
Proposition~\ref{p:innerprod-principalvectors} is that $u_i \perp v_j$
for $i \neq j$. This leads to the following important orthogonal
decomposition for the space $U+V$ based on the principal vectors.

\begin{corollary}\longthmtitle{Orthogonal Decomposition of $U+V$ by
Principal Vectors}\label{l:VU-orthogonal-decomposition}
For all $i,j \in \until{m_2}$ with $i \neq j$, we have
$[\Span(v_i)+ \Span(u_i)] \perp [\Span(v_j)+ \Span(u_j)]$. Therefore,
the space $U + V$ admits the orthogonal decomposition
\begin{align*}
U +V = \Big( \bigoplus_{k=1}^{m_2} [\Span(v_k)+ \Span(u_k)] \Big) \\
\oplus \Big( \bigoplus_{k=m_2+1}^{m_1} \!\! \Span(u_k) \Big) .   \eqoprocend
\end{align*}
\end{corollary}

The proof of this result is a direct consequence of
Lemma~\ref{l:orthonormality-principal-vectors} and
Proposition~\ref{p:innerprod-principalvectors}. Corollary~\ref{l:VU-orthogonal-decomposition}
decomposes the finite-dimensional subspace $U+V \subseteq \Fc$ into
several orthogonal spaces of dimension one or two. However, in
general, $u_i$ and $v_i$ are not orthogonal.  The next result provides
two orthonormal bases for the subspace $\Span(v_i)+ \Span(u_i)$ and
orthogonally decomposes $u_i$ and $v_i$ with respect to these bases.

\begin{lemma}\longthmtitle{Orthogonal Decomposition of Principal
Vectors}\label{l:principal-vec-decomposition}
Let $i \in \until{m_2}$ with $\theta_i \neq 0$. Then the subspace
$\Span(v_i)+ \Span(u_i)$ is two-dimensional.  Consider the
orthonormal bases\footnote{These bases can be computed using a
Gram-Schmidt process, e.g.~\citep{SJL-AB-WG:13}.}
$\{u_i,u_i^\perp\}$, $\{v_i,v_i^\perp\}$ for the subspace. Then,
\begin{subequations}
\begin{align}
v_i= \cos(\theta_i)u_i + \gamma_i \, u_i^\perp, \label{eq:vi-decomposition} 
\\
u_i= \cos(\theta_i)v_i  + \mu_i \, v_i^\perp, \label{eq:ui-decomposition} 
\end{align}
\end{subequations}
where $| \gamma_i|= | \mu_i| = \sin(\theta_i)$.
\end{lemma}
\begin{proof}
For the first part, given $\theta_i \neq 0$, we have
$\theta_i \in (0,\frac{\pi}{2}]$ based on
Definition~\ref{def:principal-angles}. Therefore,
$\innerprod{u_i}{v_i} = \innerprod{v_i}{u_i} = \cos(\theta_i) \neq
1$. Now, consider a linear combination of the form
$\alpha u_i + \beta v_i = 0$. By taking the inner product of both
sides with $u_i$ and $v_i$, one can write
\begin{align*}
\alpha + \beta \cos(\theta_i) = 0, \quad
\alpha \cos(\theta_i) + \beta = 0.
\end{align*}
Since $\cos(\theta_i) \neq 1$, the unique solution is
$\alpha = \beta = 0$ and consequently $\{u_i,v_i\}$ are linearly
independent.

To prove~\eqref{eq:vi-decomposition}, consider the following
expansion
\begin{align}\label{eq:vi-decompostion-2}
v_i = \eta_i u_i + \gamma_i u_i^\perp,
\end{align}
with $\eta_i, \gamma_i \in \cplx$.  Taking the inner product with
$u_i$, we get
\begin{align}\label{eq:a1}
\cos(\theta_i) = \eta_i.
\end{align}
Moreover, given that $u_i \perp u_i ^\perp$, one can
use~\eqref{eq:vi-decompostion-2} and the properties of norms induced
from inner products to write
\begin{align*}
\|v_i\|^2 = |\eta_i|^2 \|u_i\|^2 + |\gamma_i|^2 \|u_i^\perp\|^2.
\end{align*}
Noting that $\|v_i\| = \|u_i\| = \|u_i^\perp\| = 1$, this equality
combined with~\eqref{eq:vi-decompostion-2} and~\eqref{eq:a1} yields
$|\gamma_i| = \sin(\theta_i)$, which
proves~\eqref{eq:vi-decomposition}. The proof
of~\eqref{eq:ui-decomposition} is analogous.
\end{proof}

\section{Invariance Proximity and Principal Angles}

This section presents the main result of the paper, which provides a
closed-form formula for the invariance proximity of a subspace under
the Koopman operator using the notion of principal angles.

\begin{theorem}\longthmtitle{Closed-Form Solution for Invariance
Proximity via Principal
Angles}\label{t:invariance-proximity-closed-form} 
Let $\Sc \subseteq \Fc$ be a finite-dimensional space and let
$\Kc \Sc$ be the image of $\Sc$ under $\Kc: \Fc \to \Fc$.  Let
$0 \leq \theta_1 \leq \cdots \leq \theta_{\dim(\Kc \Sc)} \leq
\frac{\pi}{2}$ be the principal angles between $\Sc$ and $\Kc
\Sc$. Then, invariance proximity can be expressed as
\begin{align}
\Ic_{\Kc}(\Sc) = \sin\big(\theta_{\dim(\Kc \Sc)}\big).
\end{align}
Moreover, the supremum in~\eqref{def:invariance-proximity} is
actually a maximum, i.e., there exists a function $f^* \in \Sc$ such
that
\begin{align*}
\Ic_{\Kc}(\Sc) = \frac{\| \Kc f^* - \Pc_{\Sc} \Kc f^* \|}{\| \Kc f^*
\|}. 
\end{align*}
\end{theorem}
\begin{proof}
To use the results in Section~\ref{sec:prinicipal-angles}, we rely
on the following notation throughout the proof
\begin{align*}
U = \Sc,  \quad \dim(U) = m_1, \quad V = \Kc \Sc, \quad \dim(V) = m_2.
\end{align*}
Note that $m_1=\dim(\Sc) \geq \dim(\Kc \Sc)=m_2$ which is consistent
with the convention in Section~\ref{sec:prinicipal-angles}. Let then
$\{ u_1, \ldots u_{m_1} \} \subset U$ and
$\{ v_1, \ldots v_{m_2} \} \subset V$ be orthonormal bases of
principal vectors.  For convenience, for $ f \in \Sc$ with
$ \|\Kc f \| \neq 0$, we use the shorthand notation
\begin{align*}
E_\Kc(f) = \frac{\|\Kc f - \Pc_{\Sc}\Kc f \|}{\|\Kc f\|} .
\end{align*}
Our first goal is to show that $E_\Kc(f) \leq
\sin(\theta_{m_2})$. To achieve this, we decompose $\Kc f \in V$ as
\begin{align}\label{eq:Kf-decompose}
\Kc f = \sum_{i=1}^{m_2} \alpha_i v_i.
\end{align}
Since $\Pc_{\Sc}$ is the orthogonal projection on $\Sc$, we
use~\eqref{eq:Kf-decompose} in conjunction with
Proposition~\ref{p:innerprod-principalvectors} to decompose
$\Pc_{\Sc} \Kc f$ as
\begin{align}\label{eq:PKf-decompose}
\Pc_{\Sc} \Kc f
= \sum_{j =1}^{m_1} \innerprod{\Kc f}{u_j} u_j
&
=
\sum_{j=1}^{m_1} \innerprod{\sum_{i=1}^{m_2} \alpha_i v_i}{u_j}
u_j 
\nonumber
\\ 
&
= \sum_{i=1}^{m_2} \alpha_i \cos(\theta_i) u_i.
\end{align}
Using~\eqref{eq:Kf-decompose}-\eqref{eq:PKf-decompose}, the
orthogonality of subspaces in
Corollary~\ref{l:VU-orthogonal-decomposition}, and the fact that the
norm is induced by an inner product, one can write
\begin{align*}
\| \Kc f \!-\! \Pc_{\Sc} \Kc f \|^2
& \!=\! \| \sum_{i=1}^{m_2} \alpha_i
(v_i \!-\! \cos(\theta_i) u_i) \|^2  
\\
&
\!=\! \sum_{i=1}^{m_2} |\alpha_i|^2 \| v_i \!-\! \cos(\theta_i) u_i \|^2.
\end{align*}
Using now~\eqref{eq:vi-decomposition} in
Lemma~\ref{l:principal-vec-decomposition}, one can write
\begin{align}\label{eq:projection-error-closed-form}
\| \Kc f - \Pc_{\Sc} \Kc f \|^2
= \sum_{i=1}^{m_2} |\alpha_i|^2 \sin(\theta_i)^2,
\end{align}
where we have used that $\|u_i^\perp\| = 1$ for $i \in \until{m_2}$.
We can also use~\eqref{eq:Kf-decompose}, the properties of norms
induced by inner products, and $\|v_i\| = 1$ for $i \in \until{m_2}$
to write
\begin{align}\label{eq:denom-closed-form}
\| \Kc f \|^2 = \| \sum_{i=1}^{m_2} \alpha_i v_i \|^2 =
\sum_{i=1}^{m_2} |\alpha_i|^2 \|v_i\|^2 =  \sum_{i=1}^{m_2}
|\alpha_i|^2. 
\end{align}
Based
on~\eqref{eq:projection-error-closed-form}-\eqref{eq:denom-closed-form},
we have
\begin{align}\label{eq:expansion-relative-error}
(E_{\Kc}(f))^2 =
\frac{\sum_{i=1}^{m_2} |\alpha_i|^2
\sin(\theta_i)^2}{\sum_{i=1}^{m_2} |\alpha_i|^2}.
\end{align}
Since $\|\Kc f \| \neq 0$, based on~\eqref{eq:Kf-decompose}, we have
that ${\sum_{i=1}^{m_2} |\alpha_i|^2} \neq 0$. Now, since
$0 \leq \theta_1 \leq \cdots\leq \theta_{m_2} \leq \frac{\pi}{2}$,
we can write $\sin(\theta_i)^2 \leq \sin(\theta_{m_2})^2$ for all
$i \in \until{m_2}$. Therefore, for all $f \in \Sc$ with
$\|\Kc f \| \neq 0$,
\begin{align}\label{eq:relative-error-upper-bound}
(E_{\Kc}(f))^2
&=
\frac{\sum_{i=1}^{m_2} |\alpha_i|^2
\sin(\theta_i)^2}{\sum_{i=1}^{m_2} |\alpha_i|^2}
\nonumber
\\
&\leq \frac{\sum_{i=1}^{m_2} |\alpha_i|^2
\sin(\theta_{m_2})^2}{\sum_{i=1}^{m_2} |\alpha_i|^2} =
\sin(\theta_{m_2})^2  .
\end{align}
Next, we prove that the equality
in~\eqref{eq:relative-error-upper-bound} holds for some function in
$\Sc$.  Let $f^*$ belong to $ \Kc^{-1}(v_{m_2})$, the inverse image
of $v_{m_2} \in V = \Kc \Sc$ under $\Kc$  ($f^*$ exists since
$\Kc \Sc$ is the image of $\Sc$ under $\Kc$, but is generally not
unique).  Now, using~\eqref{eq:Kf-decompose} for
$\Kc f^* = v_{m_2}$, we write
$ \Kc f^* = \sum_{i=1}^{m_2} \alpha_i^* v_i$,
where $\alpha_1^* = \cdots = \alpha_{m_2-1}^*=0$ and
$\alpha_{m_2}^* = 1$. Hence, by
applying~\eqref{eq:expansion-relative-error} on $f^*$, we have
\begin{align*}
(E_{\Kc}(f^*))^2 =
\frac{\sum_{i=1}^{m_2} |\alpha_i^*|^2
\sin(\theta_i)^2}{\sum_{i=1}^{m_2} |\alpha_i^*|^2} =
\sin(\theta_{m_2})^2,
\end{align*}
and the result follows from the
definition~\eqref{def:invariance-proximity}.
\end{proof}

\begin{figure}[htb]
\centering
{\includegraphics[width=.9\linewidth]{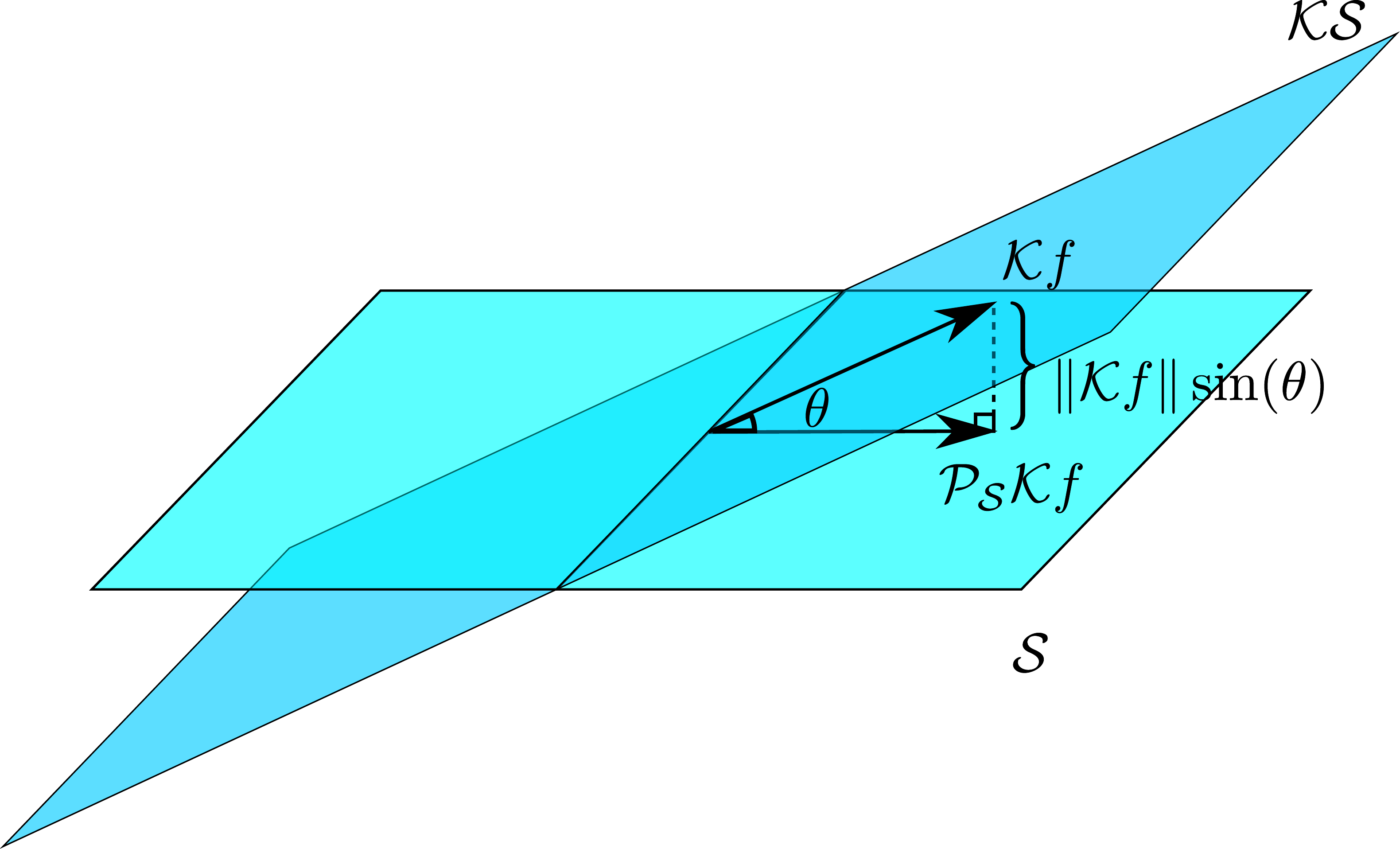}}
\caption{The error induced by the orthogonal projection on subspace
$\Sc$ is proportional to the sine of the angle between $\Kc f
\in \Kc \Sc$ and $\Pc_{\Sc} \Kc f \in \Sc$. According to
Theorem~\ref{t:invariance-proximity-closed-form}, the relative
error reaches its maximum when the angle is equal to the largest
principal angle between $\Sc$ and $\Kc
\Sc$.}\label{fig:error-angles} 
\vspace*{-1.5ex}
\end{figure}

We offer the following geometric interpretation of
Theorem~\ref{t:invariance-proximity-closed-form}.  For each function
$f \in \Sc$, note that the projection error for approximating
$\Kc f$ satisfies
$\| \Kc f- \Pc_{\Sc} \Kc f\| = \sin(\theta) \|\Kc f\|$, where
$\theta$ is the angle between $\Kc f$ and $\Pc_{\Sc} \Kc f$
(cf. Fig.~\ref{fig:error-angles}).
Theorem~\ref{t:invariance-proximity-closed-form} states that the
maximum relative error among all the functions $f$ in $\Sc$ is
achieved when the angle between $\Kc f$ and $\Pc_{\Sc} \Kc f$ is
equal to the maximum Jordan principal angle between $\Sc$ and its
image of the Koopman operator $\Kc \Sc$.

\begin{remark}\longthmtitle{Invariance Proximity is not a Metric} {\rm
Note there is a difference between invariance proximity and
notions of metrics for subspaces based on principal angles.
Theorem~\ref{t:invariance-proximity-closed-form} might create the
illusion that invariance proximity is a gap metric, see
e.g.~\citep{MIK:1975,AK-TTG:21}. However, this is not generally
true since the dimensions of $\Sc$ and $\Kc \Sc$ might be
different and we \emph{only} project from $\Kc \Sc$ onto $\Sc$.
\oprocend }
\end{remark}

Theorem~\ref{t:invariance-proximity-closed-form} provides a
closed-form expression for invariance proximity based on the
well-known concept of principal angles and vectors, providing a direct
insight into the geometry of projection-based Koopman
methods. Moreover, the algebraic decomposition based on principal
vectors paves the way for the direct calculation of invariance
proximity, which is what we discuss in the next section.

\section{Numerical Computation of Invariance Proximity}
The computation of principal angles for invariance proximity relies on
the optimization on function spaces, which requires calculation of
multi-variable integrals. Here, we provide a transformation to compute
invariance proximity using efficient numerical linear algebra.  To
achieve this, we first embed all the subspaces of interest into a
larger finite-dimensional subspace of $\Fc$.

\begin{lemma}\longthmtitle{Finite-dimensional Subspace
Embedding}\label{l:finite-dimensional-embedding}
Let $\Sc \subset \Fc$ be a finite-dimensional subspace and define
$\Wc = \Sc + \Kc \Sc$. Then, $\Wc$ is finite-dimensional and complete
(in the metric induced by the inner product).
\end{lemma}
\begin{proof}
Since both $\Sc$ and $\Kc \Sc$ are finite dimensional, their sum is
also finite-dimensional. The second part directly
follows~\citep[Theorem~2.4-2]{EK:89}.
\end{proof}

Lemma~\ref{l:finite-dimensional-embedding} shows that $\Wc$ in its own
right is a Hilbert space. Noting that $\Sc, \Kc \Sc \subset \Wc$, we
connect $\Wc$ to a more suitable subspace for numerical computations
via an isomorphism\footnote{An isomorphism between two Hilbert spaces
is a linear bijective map that preserves the inner product (and
induced norm and metric)~\citep[Section~3.2]{EK:89}.}.

\begin{lemma}\longthmtitle{Isomorphism between $\Wc$ and
$\cplx^{\dim(\Wc)}$}\label{l:isomorphism}
Consider\footnote{If the function space $\Fc$ is defined over field
$\real$, one can similarly build an isomorphism between $\Wc$ and
$\real^{\dim(\Wc)}$ and all the ensuing results will remain valid
given this change.} the space $\cplx^{\dim(\Wc)}$ endowed with the
standard inner product
$\innerprod{\cdot}{\cdot}_{\cplx^{\dim(\Wc)}}$.  Let
$\{w_1, \ldots, w_{\cplx^{\dim(\Wc)}}\}$ and
$\{c_1, \ldots, c_{\cplx^{\dim(\Wc)}}\}$ be \emph{orthonormal} bases
for $\Wc$ and $\cplx^{\dim(\Wc)}$, respectively. Define the
\emph{linear} map $Q: \Wc \to \cplx^{\dim(\Wc)}$ such that
$w_i \mapsto c_i$ for all $i \in \until{\dim(\Wc)}$. Then, $Q$ is an
isomorphism, i.e.,
\begin{enumerate}
\item $Q$ is bijective;
\item $\innerprod{x}{y} = \innerprod{Qx}{Qy}_{\cplx^{\dim(\Wc)}}$,
$\forall x,y \in \Wc$;
\item
$\innerprod{m}{n}_{\cplx^{\dim(\Wc)}} = \innerprod{Q^{-1}
m}{Q^{-1} n}$, $\forall m,n \in \cplx^{\dim(\Wc)}$. \oprocend
\end{enumerate}
\end{lemma}

The construction of $Q$ in Lemma~\ref{l:isomorphism} is well-known in
the literature (see e.g., the proof
of~\citep[Theorem~3.6-5(b)]{EK:89}) and the proof is a direct
consequence of linearity of $Q$ and the properties of inner products.

Based on Lemma~\ref{l:isomorphism}, the spaces $\Wc$ and
$\cplx^{\dim(\Wc)}$ have the same structure and any algorithmic
operation involving inner products (and induced norms and metrics) has
the same effect in both spaces.  Since the definition of principal
angles only depends on inner products and induced norms, they are
preserved under the isomorphism.

\begin{corollary}\longthmtitle{Isomorphisms Preserve Principal
Angles}\label{c:iso-principal-angle}
Let $U,V \subset \Wc$ and $Q(U)$ and $Q(V)$ be the images of $U$ and
$V$ under the isomorphism $Q: \Wc \to \cplx^{\dim(\Wc)}$. Then, the
principal angles between $U, V \subset \Wc$ and
$Q(U),Q(V) \subset \cplx^{\dim(\Wc)}$ are the same. \oprocend
\end{corollary}

Corollary~\ref{c:iso-principal-angle} has important practical
consequences since it allows one to compute the principal angles
between subspaces of general inner-product spaces using efficient
numerical algorithms developed for the special case of
$\cplx^{\dim(\Wc)}$. This allows a universal formulation without the
need to design specific algorithms based the choice of inner product
space. We next use Corollary~\ref{c:iso-principal-angle} to provide a
revised version of Theorem~\ref{t:invariance-proximity-closed-form}
that relies on computations in $\cplx^{\dim(\Wc)}$.

\begin{theorem}\longthmtitle{Invariance Proximity via
Isomorphisms}\label{t:invariance-proximity-isomorphism}
Let $\Sc \subseteq \Fc$ be a finite-dimensional space and $\Kc \Sc$
be its image under $\Kc: \Fc \to \Fc$.  Let $\Wc= \Sc + \Kc \Sc$ and
consider the isomorphism $Q: \Wc \to \cplx^{\dim(\Wc)}$. Also, let
$0 \leq \gamma_1 \leq \cdots \leq \gamma_{\dim(Q(\Kc \Sc))} \leq
\frac{\pi}{2}$ be the principal angles between $Q(\Sc)$ and
$Q(\Kc \Sc)$. Then,
$ \Ic_{\Kc}(\Sc) = \sin\big(\gamma_{\dim(Q(\Kc
\Sc))}\big)$. \oprocend
\end{theorem}

\begin{remark}\longthmtitle{Numerical Computation of Invariance
Proximity}\label{r:numerical-computation}
{\rm Based on Theorem~\ref{t:invariance-proximity-isomorphism},
one can compute invariance proximity by finding the principal
angles between subspaces comprised on $n$-tuples of numbers,
instead of directly working with functions. There exist
efficient routines for this purpose, e.g.~\citep{AB-GHG:73},
which require finding orthonormal basis for subspaces $Q(\Sc)$
and $Q(\Kc \Sc)$ in
Theorem~\ref{t:invariance-proximity-isomorphism}, and finding
the maximum singular value of a matrix with the dimensions of
order $\dim(Q(\Sc))$. This can be done through truncated
Singular Value Decomposition (SVD) on matrices with time
complexity of $O(\dim(\Sc)^3)$ FLOPs
(see,~e.g.~\citep{XL-SW-YC:19}). It is worth mentioning that
MATLAB\textsuperscript{\textregistered} has a built-in command
based on the algorithm in~\citep{AB-GHG:73}.  In case the
principal angles are small, this algorithm can struggle due to
round-off errors. The work in~\citep{AVK-MEA:02} provides a
modified algorithm to address this issue, which has been
implemented in the SciPy package.  \oprocend }
\end{remark}

Next, we provide Algorithm~\ref{algo:invariance_proximity} to
encapsulate how to apply the paper's results in order to numerically
compute invariance proximity in general inner-product spaces. Even
though our treatment considers complex inner-product spaces, all the
results identically apply to the (simpler) case of real
inner-product spaces (one only need to replace the complex
inner-product with a real one).  For this reason,
Algorithm~\ref{algo:invariance_proximity} considers inner-product
spaces on field~$\Fbb$ (which can be either $\cplx$ or $\real$).

\begin{algorithm}[!htb]
\caption{Computing Invariance Proximity} \label{algo:invariance_proximity}
\begin{algorithmic}[1] 
\Statex \textbf{Inputs:} 
\Statex $\bullet$ Basis $\Psi = [\Psi_1, \ldots, \Psi_{\dim(\Sc)}]$ for subspace $\Sc$ 
\Statex $\bullet$ Map $T$ for system~\eqref{eq:dynamical-sys} or operator $\Kc$ in~\eqref{eq:Koopman-def}
\Statex $\bullet$  Field of the function space denoted by $\Fbb$ (either $\cplx$ or $\real$)
\Statex \textbf{Output:} $\Ic_{\Kc}(\Sc)$ \Comment{Invariance proximity of $\Sc$}
\smallskip \smallskip
\Statex \textbf{Procedure:}
\smallskip
\Statex $\triangledown$ Find a generator for $\Kc \Sc$ (might not be linearly independent)
\State $Z = [\Kc \Psi_1, \ldots, \Kc \Psi_{\dim(\Sc)}]$
\medskip
\Statex $\triangledown$ Find an orthonormal basis for $\Kc \Sc$ 
\State $\Phi = [\Phi_1, \ldots,  \Phi_{\dim(\Kc \Sc)}]:= \operatorname{orth}(Z)$
\medskip
\Statex $\triangledown$ Concatenate $\Psi$ and $\Phi$
\State $[\Psi, \Phi] =  [\Psi_1, \dots, \Psi_{\dim(\Sc)}, \Phi_1, \ldots, \Phi_{\dim(\Kc \Sc)}]$ 
\medskip
\Statex $\triangledown$ Find an orthonormal basis for $\Wc = \Sc + \Kc \Sc $
\State $W = [w_1, \ldots, w_{\dim(\Wc)}]:=\operatorname{orth}([\Psi, \Phi]) $ 
\medskip
\Statex $\triangledown$ Form isomorphism between $\Wc$ and $\Fbb^{\dim(\Wc)}$ (cf.~Lemma~\ref{l:isomorphism})
\State  Define $Q$ by $\sum_{i=1}^{\dim(\Wc)} \alpha_i w_i \mapsto \sum_{i=1}^{\dim(\Wc)} \alpha_i e_i$  \Comment $\{e_1, \ldots e_{\dim(\Wc)} \}$ is an orthonormal basis for $\Fbb^{\dim(\Wc)}$
\medskip
\Statex $\triangledown$ Form matrix $Q_{\Sc}^\text{basis} \in \Fbb^{\dim(\Wc) \times \dim(\Sc)}$ whose columns span $Q(\Sc)$
\State $Q_{\Sc}^\text{basis} = [Q(\Psi_1), \dots, Q(\Psi_{\dim(\Sc)})]$ 
\medskip
\Statex $\triangledown$ Form matrix $Q_{\Kc \Sc}^\text{basis} \in \Fbb^{\dim(\Wc) \times \dim(\Kc \Sc)}$ whose columns span $Q(\Kc \Sc)$
\State $Q_{\Kc \Sc}^\text{basis} = [Q(\Phi_1), \dots, Q(\Phi_{\dim(\Kc \Sc)})]$
\medskip
\Statex $\triangledown$ Compute maximum principal angle between $Q(\Sc)$ and $Q(\Kc \Sc)$ by applying numerical routines on $Q_{\Sc}^\text{basis}$ and $Q_{\Kc \Sc}^\text{basis}$ (cf.~Remark~\ref{r:numerical-computation})
\State $\gamma_{\dim(Q(\Kc \Sc))} = \text{maximum principal angle}$
\medskip
\Statex $\triangledown$ Compute invariance proximity (Theorem~\ref{t:invariance-proximity-isomorphism})
\State $\Ic_{\Kc}(\Sc) = \sin(\gamma_{\dim(Q(\Kc \Sc))})$
\end{algorithmic}
\end{algorithm}

In Algorithm~\ref{algo:invariance_proximity}, the function
$\texttt{orth}$ provides an orthonormal basis for its argument, which
can be computed via a Gram-Schmidt process
(e.g.~\citep{SJL-AB-WG:13}) and removing the redundant (linearly
dependent) terms.

We finish the section by observing that Theorem~6.4 can also be used
for data-driven cases.

\begin{remark}\longthmtitle{Computing Invariance Proximity for
Data-Driven Cases} {\rm In data-driven settings, generally the
system and its associated Koopman operator are unknown and only
data snapshots of trajectories are available. In such cases, the
inner product is empirically defined based on the data, e.g., the
well-known EDMD method~\citep{MOW-IGK-CWR:15,MK-IM:18}. Given that
such inner products only depend on the value of functions over a
data set, to apply
Theorem~\ref{t:invariance-proximity-isomorphism}, one does not
need full knowledge of the elements in $\Kc \Sc$: instead, their
value on the data set is enough. Such values can be computed
since, for $ \Kc g \in \Kc \Sc$ with $g \in \Sc$, we have
$\Kc g (x) = g \circ T(x) = g(x^+)$, where $x^+= T(x)$ is the next
point on the trajectory from~$x$. \oprocend }
\end{remark}

\section{Simulation Results}
Let the system with state $x = [x_1,x_2]^T$ on $\Xc = [-1,1]^2$,
\begin{align}\label{eq:nonlinear-sys-example}
x_1^+ &= 0.9 \, x_1 
\nonumber \\
x_2^+ &= 0.4 \, \big(\sin(x_2) + x_1^2 \big) + 0.01 \, x_2^2.
\end{align}
Consider the function space $\Fc$ (over $\real$) comprised of all
real-valued continuous functions with domain $\Xc$, equipped with the
inner product $\innerprod{f}{g} = \int \int_{\Xc} f(x) g(x) dx_1 dx_2$
for $f,g \in \Fc$. We compute the invariance proximity for
subspaces\footnote{Note that the elements in the set are
functions. For example $x_1$ represents $f(x)=x_1$. This is a
conventional notation in the literature.}
$\Sc_1 = \Span\{1, x_1, x_1^2\}$,
$\Sc_2 = \Span\{1, x_1, x_2, x_1^2\}$, and
$\Sc_3 = \Span\{1, x_1, x_2, x_1^2, x_2^2\}$. We explain the procedure
for building the model and finding the invariance proximity for
subspace $\Sc_1$. The procedure for other subspaces is identical.

\emph{Finding the model:} (i) we first apply the Gram-Schmidt process
on the basis of $\Sc_1$ and create a function $\Psi: \Xc \to \real^3$
whose elements form an orthonormal basis for $\Sc_1$; (ii) we apply
the Koopman operator on $\Psi$
following~\eqref{eq:restriction-on-noninvariant-basis} and find the
matrix $\Kapprox$ whose $ij$th element can be computed by
$[\Kapprox]_{ij} = \innerprod{\Kc \Psi_i}{\Psi_j}$, where $\Psi_i$ and
$\Psi_j$ are the $i$th and $j$th elements of $\Psi$, resp.

\emph{Computing the invariance proximity:} To compute invariance
proximity, we employ
Algorithm~\ref{algo:invariance_proximity}. Here, we discuss the
implementation details: (i) we find a basis for space $\Kc \Sc_1$
by applying $\Kc$ on the elements of a basis for $\Sc_1$, then
performing a Gram-Schmidt process and removing the redundant (linearly
dependent) terms;
(ii) we find an orthonormal basis for $\Wc = \Sc_1 + \Kc \Sc_1$ by
concatenating the basis elements of $\Sc_1$ and $\Kc \Sc_1$, and
applying the Gram-Schmidt algorithm and removing the linearly
dependent elements. We denote this basis by
$\{w_1, \ldots, w_{\dim(\Wc)}\}$;\footnote{In
Steps~(i)-(ii), we relied on the symbolic computations in
MATLAB\textsuperscript{\textregistered}. Alternatively, one can
compute the inner-product through numerical integration.}
(iii) we define the isomorphism $Q: \Wc \to \real^{\dim (\Wc)}$
(cf.~Lemma~\ref{l:isomorphism}) by mapping
$\sum_{i=1}^{\dim(\Wc)} \alpha_i w_i \mapsto \sum_{i=1}^{\dim(\Wc)}
\alpha_i e_i = [\alpha_1, \ldots, \alpha_{\dim(\Wc)}]^T$, where $e_i$
is the $i$th element of canonical basis ($i$th column of the identity
matrix $I_{\dim(\Wc)}$) for $\real^{\dim(\Wc)}$; (iv) to apply
Theorem~\ref{t:invariance-proximity-isomorphism}, we find bases for
subspaces $Q(S_1)$ and $Q(\Kc \Sc_1)$. Given the basis $\Psi$ of
$\Sc_1$, we compute the action of $Q$ on elements of $\Psi$ by
$[\innerprod{\Psi_i}{w_1}, \ldots,
\innerprod{\Psi_i}{w_{\dim(\Wc)}}]^T$.  We concatenate these vectors
into a matrix
$Q_{\Sc_1}^{\text{basis}} \in \real^{\dim(\Wc) \times \dim(\Sc_1)}$,
whose range space is $Q(S_1)$. Similarly, we form the matrix
$Q_{\Kc \Sc_1}^{\text{basis}} \in \real^{\dim(\Wc) \times \dim(\Kc
\Sc_1)}$, whose range space is $Q(\Kc \Sc_1)$; (v) finally, to
invoke Theorem~\ref{t:invariance-proximity-isomorphism}, we use the
built-in \texttt{subspace} command in
MATLAB\textsuperscript{\textregistered} (which is based
on~\citep{AB-GHG:73}) to compute the maximum principle angle
between range spaces of $Q_{\Sc_1}^{\text{basis}}$ and
$Q_{\Kc \Sc_1}^{\text{basis}}$. By
Theorem~\ref{t:invariance-proximity-isomorphism}, the invariance
proximity equals the sine of this angle.

\emph{Interpretation and discussion:}
Table~\ref{table:invariance-proximity} shows the invariance proximity
for subspaces $\Sc_1$, $\Sc_2$, and $\Sc_3$. Clearly, $\Sc_1$ is
Koopman invariant since its functions are monomials of the first state
variable $x_1$ and the evolution of $x_1$ abides by a linear dynamics
($x_1^+ = 0.9 \, x_1$). This is consistent with the fact that
invariance proximity of $\Sc_1$ is zero\footnote{Given that the
computation is done by a digital computer, the invariance
proximity is at the level of machine precision instead of exact
zero.}. The invariance proximity for $\Sc_2$ is rather small,
indicating that the worst-case relative function prediction error
(note that Koopman operator acts on functions) is $4.8 \%$. On the
other hand, the worst-case relative function prediction error for
$\Sc_3$ is $82.3 \%$, rendering models on $\Sc_3$ unreliable. This is
despite the fact that $\Sc_2 \subset \Sc_3$, which indicates that a
larger subspace is not necessarily better\footnote{This does not
contradict the asymptotic results in~\citep{MK-IM:18}. See
\citep[Example~2.1]{MH-JC:23-auto} for a discussion.}.

{ 
\renewcommand{\arraystretch}{1.5}
\begin{table}[htb]
\centering
\caption{Invariance proximity for subspaces $\Sc_1$, $\Sc_2$, and $\Sc_3$.} \label{table:invariance-proximity}
\begin{tabular}[htb]{ | l || c | c | c |  }
\hline
\textbf{Subspace}	               			  & $\Sc_1$  & $\Sc_2$  & $\Sc_3$   \\ \hline %
\textbf{Invariance Proximity}           & $\sim 0$       & 0.048     &    0.823        \\ \hline %
\end{tabular}
\vspace*{-1ex}
\end{table}
} 

\begin{figure}[htb]
\centering 
{\includegraphics[width=.48\linewidth]{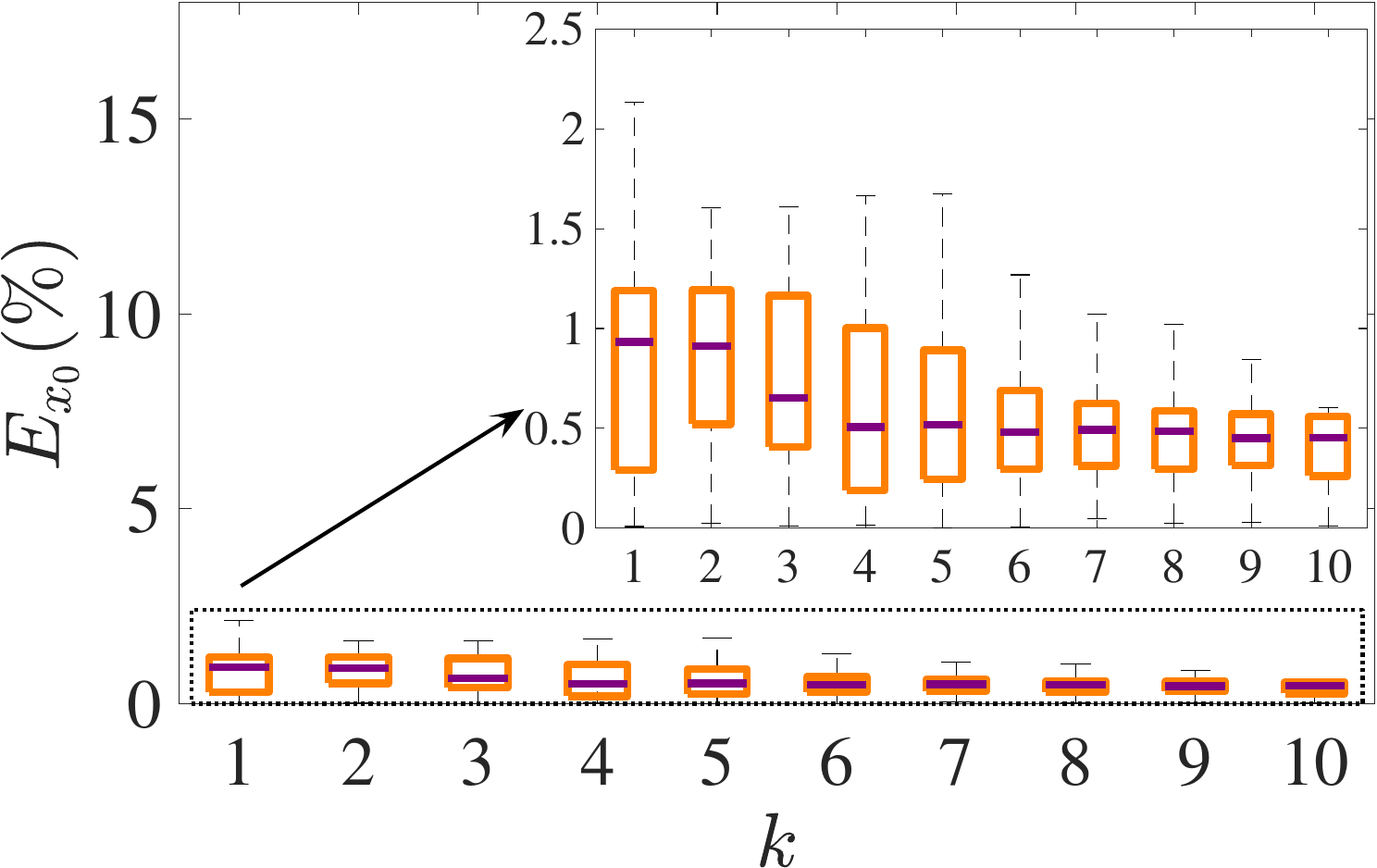}}
{\includegraphics[width=.48\linewidth]{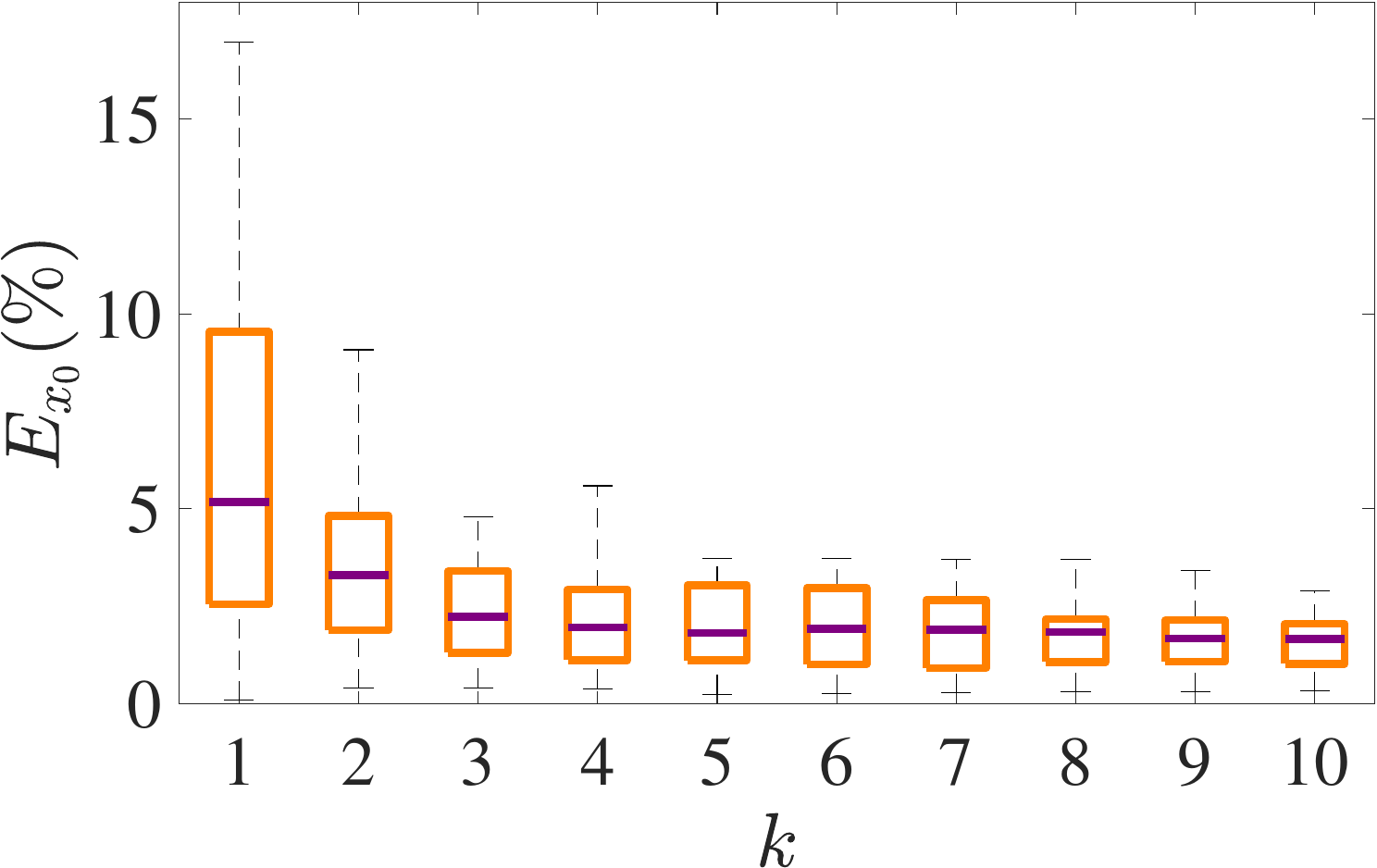}}
\caption{The median (in purple), the distance between 25\% and 75\%
percentiles (orange box), and the entire range (whiskers) for the
relative error in~\eqref{eq:relative-error} over 100 trajectories
with initial conditions uniformly sampled from $\Xc$ given the
projected models for subspaces $\Sc_2$ (left) and $\Sc_3$
(right).}\label{fig:relative-error}
\vspace*{-1.5ex}
\end{figure}

To show how the subspace's quality impacts the accuracy of the linear
predictor in~\eqref{eq:restriction-on-noninvariant-basis} on the
system trajectories, we consider the following relative error function
given a trajectory $\{x(k)\}_{k \in \naturals_0}$ from the initial
condition $x_0$
\begin{align}\label{eq:relative-error}
E_{x_0}(k) = \frac{\| \Psi(x(k)) - \Kapprox^k \Psi(x_0) \|}{\|
\Psi(x(k)) \|} \times 100. 
\end{align}
In the equation above, $\| \cdot \|$ denotes the usual Euclidean
norm in $\real^n$. Unlike invariance proximity, which does not
depend on the choice of basis, the error
in~\eqref{eq:relative-error} is influenced by the basis $\Psi$ when
it is not orthonormal. Therefore, to make a fair portrayal, we
enforce the elements of $\Psi$ to be orthonormal. The error
in~\eqref{eq:relative-error} for subspace $\Sc_1$ is equal to zero for
all initial conditions since the subspace is Koopman invariant and the
prediction is exact. Fig.~\ref{fig:relative-error} shows the error
in~\eqref{eq:relative-error} for subspaces $\Sc_2$ and $\Sc_3$ over
100 system trajectories with the length of $10$ time steps and initial
conditions uniformly sampled from $\Xc =
[-1,1]^2$. Fig.~\ref{fig:relative-error} clearly shows the superiority
of the model on $\Sc_2$ compared to $\Sc_3$. Also, it is worth
mentioning that the variance of the error is much lower for subspace
$\Sc_2$ since invariance proximity indicates the function prediction
errors over the entire state space instead of a single or a few
initial conditions.

\section{Conclusions}
We have provided a closed-from description of invariance proximity, a
notion that measures the worst-case relative error of Koopman-based
projected models, over general complex inner-product spaces.  Our
solution leverages the geometry behind projections, subspaces, and the
Koopman operator and relies on the calculation of Jordan principal
angles between two finite-dimensional subspaces in an
infinite-dimensional space of functions.  To avoid the computational
challenges of calculating inner products and principal angles in such
a space, we have used specific isomorphisms to make the problem of
computing invariance proximity amenable to efficient algorithmic
routines from numerical linear algebra. Future work will include using
invariance proximity to provide (a) performance guarantees on
Koopman-based prediction and control schemes, and (b) safety and
stability certificates from Koopman-based approximate models.

\section*{Acknowledgments}
This work was supported by ONR Award N00014-23-1-2353.  The
authors would like to thank the anonymous reviewers for their
insightful comments and suggestions, which helped enhance the
quality of the paper.

\end{document}